\documentclass{article}

\usepackage{lmodern}
\usepackage{amsthm,amsmath,amssymb,oldgerm,hyperref,mathrsfs,amsxtra, bm, accents}
\usepackage[toc,page]{appendix}
\usepackage[a4paper, margin=2.5 cm]{geometry}
\usepackage{mathtools}
\usepackage{dsfont}
\usepackage{bbm} 
\usepackage{theoremref}
\usepackage{amsopn}
\usepackage{aliascnt} 
\usepackage{enumerate}
\usepackage{stmaryrd} 
\usepackage{afterpage} 
\usepackage{extarrows} 
\usepackage{commath} 
\usepackage{pgf,tikz-cd}
\usepackage{mathrsfs}
\usetikzlibrary{arrows}
\usepackage{esint} 
\newcommand{\mynewtheorem}[2]{
\newaliascnt{#1}{dummy}
\newtheorem{#1}[#1]{#2}
\aliascntresetthe{#1}
\expandafter\def\csname #1autorefname\endcsname{#2}}
 
\theoremstyle{plain}
\mynewtheorem{theorem}{Theorem}
\mynewtheorem{proposition}{Proposition}
\mynewtheorem{corollary}{Corollary}
\mynewtheorem{lemma}{Lemma}

\newtheorem*{theorem*}{Theorem}
\theoremstyle{definition}
\mynewtheorem{definition}{Definition}
\mynewtheorem{example}{Example}
\theoremstyle{remark}
\mynewtheorem{remark}{Remark}

\def\equationautorefname~#1\null{equation~(#1)\null}

\numberwithin{equation}{section}
\DeclareMathOperator{\im}{Im}

\DeclareMathOperator{\tr}{\mathrm{tr}}

\begin{document}

\title{Relating Siegel cusp forms to Siegel--Maa\ss{} forms}
\author{J\"urg Kramer and Antareep Mandal\footnote{The first author acknowledges support from the DFG Cluster of Excellence MATH+.}}
\date{\today}

\maketitle

\begin{abstract}
\noindent
In this paper we generalize a well-known isomorphism between the space of cusp forms of weight $k$ for a Fuchsian subgroup of the first kind $\Gamma
\subset\mathrm{SL}_{2}(\mathbb{R})$ and the space of certain Maa\ss{} forms of weight $k$ for $\Gamma$ to an isomorphism between the space of 
Siegel cusp forms of weight $k$ for a subgroup $\Gamma\subset\mathrm{Sp}_{n}(\mathbb{R})$, which is commensurable with the Siegel modular group 
$\mathrm{Sp}_{n}(\mathbb{Z})$, and a suitable space of Siegel--Maa\ss{} forms of weight $k$ for $\Gamma$.
\end{abstract}

\section{Introduction}
Let $\mathbb{H}:=\{z=x+i y\in\mathbb{C}\,\vert\,y>0\}$ denote the upper half-plane and $\Gamma\subset\mathrm{SL}_{2}(\mathbb{R})$ a Fuchsian 
subgroup of the first kind, which acts by fractional linear transformations on $\mathbb{H}$. Let $\mathcal{S}_{k}(\Gamma)$ denote the space of cusp  
forms of weight $k$ for $\Gamma$ and let $\mathcal{H}_{k}(\Gamma)$ denote the space of real-analytic automorphic forms of weight $k$ for $\Gamma$, 
on which the Maa\ss{} Laplacian 
\begin{align*}
\Delta_{k}:=y^{2}\bigg(\frac{\partial^{2}}{\partial x^{2}}+\frac{\partial^{2}}{\partial y^{2}}\bigg)-iky\frac{\partial}{\partial x}
\end{align*}
of weight $k$ acts. Then, it is well-known that there is an isomorphism
\begin{align}
\label{1-iso}
\mathcal{S}_{k}(\Gamma)\cong\ker\bigg(\Delta_{k}+\frac{k}{2}\bigg(1-\frac{k}{2}\bigg)\mathrm{id}\bigg)
\end{align}
of $\mathbb{C}$-vector spaces, induced by the assignment $f\mapsto y^{k/2}f$, where the right-hand side consists of the Maa\ss{} forms in $\mathcal
{H}_{k}(\Gamma)$ with eigenvalue $k/2(1-k/2)$ of $\Delta_{k}$. This identification of two types of automorphic forms for $\Gamma$ has various useful 
applications. For example, in the article \cite{JK2}, the isomorphism~\eqref{1-iso} was crucial in relating the sup-norm bound problem for cusp forms of 
weight $k$ for $\Gamma$ to bounds for the heat kernel for $\Delta_{k}$ on the quotient space $\Gamma\backslash\mathbb{H}$. 

In this paper, we attempt a generalization of the isomorphism~\eqref{1-iso} to the Siegel modular setting, which, to our surprise, we could not find in 
the literature. Letting $\mathrm{Sym}_{n}(\mathbb{C})$ be the set of complex symmetric $(n\times n)$-matrices, we let $\mathbb{H}_{n}:=\{Z=X+iY\in
\mathrm{Sym}_{n}(\mathbb{C})\,\vert\,Y>0\}$ denote the Siegel upper half-space of degree $n\geq 1$ and we let $\Gamma\subset\mathrm{Sp}_{n}
(\mathbb{R})$ denote a subgroup acting by generalized fractional linear transformations on $\mathbb{H}_{n}$, which is commensurable with the Siegel 
modular group $\mathrm{Sp}_{n}(\mathbb{Z})$. Then, let $\mathcal{S}_{k}^{n}(\Gamma)$ denote the space of Siegel cusp forms of weight $k$ and 
degree $n$ for $\Gamma$ and let $\mathcal{H}_{k}^{n}(\Gamma)$ be the space of real-analytic automorphic forms of weight $k$ and degree $n$ for 
$\Gamma$, on which the Siegel--Maa\ss{} Laplacian
\begin{align}
\label{deltak}
\Delta_{k}:=\tr\bigg(Y\bigg(\bigg(Y\dpd{}{X}\bigg)^{t}\dpd{}{X}+\bigg(Y\dpd{}{Y}\bigg)^{t}\dpd{}{Y}\bigg)-ikY\dpd{}{X}\bigg)
\end{align}
of weight $k$ acts. As the main result of this paper, we show in Corollary~\ref{kernel connection} the isomorphism
\begin{align}
\label{1-iso general}
\mathcal{S}_{k}^{n}(\Gamma)\cong\ker\bigg(\Delta_{k}+\frac{nk}{4}(n-k+1)\mathrm{id}\bigg),
\end{align}
of $\mathbb{C}$-vector spaces, induced by the assignment $f\mapsto\det(Y)^{k/2}f$, thereby generalizing the  isomorphism~\eqref{1-iso} to the Siegel 
modular setting. The right-hand side of~\eqref{1-iso general} now consists of the Siegel--Maa\ss{} forms in $\mathcal{H}_{k}^{n}(\Gamma)$ with 
eigenvalue $nk(n-k+1)/4)$ of $\Delta_{k}$.

In case $n=1$, the isomorphism~\eqref{1-iso} is obtained as a by-product of the proof of the symmetry of the Maa\ss{} Laplacian $\Delta_{k}$ (see
\cite{Roelcke}). The most straightforward proof of the symmetry of $\Delta_{k}$ is obtained by constructing a suitable $\mathrm{SL}_{2}(\mathbb{R})
$-invariant differential form using the raising or the lowering operators on $\mathbb{H}$, and then integrating it over the quotient space $\Gamma
\backslash\mathbb{H}$ (see \cite{Bump}, p.~135). Generalizations of all these operators as well as their transformation behaviour under the action 
of the symplectic group $\mathrm{Sp}_{n}(\mathbb{R})$ to the Siegel modular setting have been provided by Maa\ss{} in~\cite{Maass2}. However, 
in spite of all these crucial ingredients being around for a long time, we could not find in the literature a precise proof of the symmetry of the 
Siegel--Maa\ss{} Laplacian $\Delta_{k}$ of weight $k$. We provide a complete proof of the symmetry of $\Delta_{k}$ in Theorem~\ref{main theorem}, 
where we construct the appropriate $\mathrm{Sp}_{n}(\mathbb{R})$-invariant differential form on $\mathbb{H}_{n}$, which, while calculationally a bit 
cumbersome, is conceptually a rather straightforward piecing-together of Maa\ss{}'s calculations. Our main result in Corollary~\ref{kernel connection} 
is then a consequence of Theorem~\ref{main theorem}.

As indicated above, the generalization of the isomorphism~\eqref{1-iso} will perspectively allow us, among others, to relate the sup-norm bound 
problem for Siegel cusp forms of weight $k$ and degree $n$ for $\Gamma$ to bounds for the heat kernel for the Siegel--Maa\ss{} Laplacian $\Delta_
{k}$ on the quotient space $\Gamma\backslash\mathbb{H}_{n}$. 

This paper is organized as follows: In section~2, we provide a quick summary of the basics of the Siegel upper half-space and Siegel modular forms. 
In the subsequent two sections~3 and~4, we introduce and discuss the transformation behaviour of the relevant operators in the Siegel modular setting. 
This material is already present in~\cite{Maass2}, but due to sub-optimal typesetting, at places, it is hard to decipher. So we take this opportunity to 
redo these calculations along Maa\ss{}'s lines and present them here for the reader's convenience. However, no claim of originality is made here on 
this material. Finally in section~5, piecing together Maa\ss{}'s results, we construct the appropriate $\mathrm{Sp}_{n}(\mathbb{R})$-invariant differential 
form on $\mathbb{H}_{n}$ to show the symmetry of the Siegel--Maa\ss{} Laplacian $\Delta_{k}$, and then use it to show the generalization~\eqref
{1-iso general} of the isomorphism~\eqref{1-iso}.

\section{Basic notations and definitions}
For $n\in\mathbb{N}_{>0}$ and a commutative ring $R$, let $\mathrm{M}_{n}(R)$ denote the set of $(n\times n)$-matrices with entries in $R$ and 
$\mathrm{Sym}_{n}(R)$ the set of symmetric matrices in $\mathrm{M}_{n}(R)$. The Siegel upper half-space 
$\mathbb{H}_{n}$ of degree $n$ is then defined by
\begin{align*}
\mathbb{H}_{n}:=\{Z=X+iY\in\mathrm{M}_{n}(\mathbb{C})\,\vert\,X,Y\in\mathrm{Sym}_{n}(\mathbb{R}):\,Y>0\}.
\end{align*}
The symplectic group $\mathrm{Sp}_{n}(\mathbb{R})$ of degree $n$ is defined by 
\begin{align*}
\mathrm{Sp}_{n}(\mathbb{R}):=\{\gamma\in\mathrm{M}_{2n}(\mathbb{R})\,\vert\,\gamma^{t}J_{n}\gamma=J_{n}\},
\end{align*}
where $J_{n}\in\mathrm{M}_{2n}(\mathbb{R})$ is the skew-symmetric matrix 
\begin{align*}
J_{n}:=\begin{pmatrix}0&\mathbbm{1}_{n}\\-\mathbbm{1}_{n}&0\end{pmatrix}
\end{align*}
with $\mathbbm{1}_{n}$ denoting the identity matrix of $\mathrm{M}_{n}(\mathbb{R})$. Writing an element $\gamma\in\mathrm{Sp}_{n}(\mathbb{R})$ 
in block form as
\begin{align*}
\gamma=\begin{pmatrix}A&B\\C&D\end{pmatrix}
\end{align*}
with $A,B,C,D\in\mathrm{M}_{n}(\mathbb{R})$, we can recast the relation $\gamma^{t}J_{n}\gamma=J_{n}$ as the set of relations 
\begin{align}
\label{symprelation 1}
A^{t}C=C^{t}A,\quad B^{t}D=D^{t}B,\quad A^{t}D-C^{t}B=\mathbbm{1}_{n}.
\end{align}
Observing that with $\gamma\in\mathrm{Sp}_{n}(\mathbb{R})$, we also have $\gamma^{t}\in\mathrm{Sp}_{n}(\mathbb{R})$, the set of relations~\eqref
{symprelation 1} turns out to be equivalent to the relations
\begin{align*}
 AB^{t}=BA^{t},\quad CD^{t}=DC^{t},\quad AD^{t}-BC^{t}=\mathbbm{1}_{n}.  
\end{align*}

The group $\mathrm{Sp}_{n}(\mathbb{R})$ acts by the symplectic action
\begin{align}
\label{sympaction}
\mathbb{H}_{n}\ni Z\mapsto \gamma Z=(AZ+B)(CZ+D)^{-1} \qquad\big(\gamma=\big(\begin{smallmatrix}A&B\\C&D\end{smallmatrix}\big)\in\mathrm{Sp}_
{n}(\mathbb{R})\big)
\end{align}
on $\mathbb{H}_{n}$. Using the equality
\begin{align}
\label{actsymmetry}
(AZ+B)(CZ+D)^{-1}=(CZ+D)^{-t}(AZ+B)^{t},
\end{align}
which follows from equation~\eqref{symprelation 1}, we compute 
\begin{align*}
(CZ+D)^{t}\,\mathrm{Im}(\gamma Z)(C\overline{Z}+D)=\mathrm{Im}(Z),
\end{align*}
from which we derive the important relation
\begin{align}
\label{imaginary transform}
\mathrm{Im}(\gamma Z)=(CZ+D)^{-t}\,\mathrm{Im}(Z)(C\overline{Z}+D)^{-1}
\end{align}
giving rise to the determinant relation
\begin{align*}
\det(\mathrm{Im}(\gamma Z))=\frac{\det(\mathrm{Im}(Z))}{\vert\det(CZ+D)\vert^{2}}.
\end{align*}

For the differential of the symplectic action~\eqref{sympaction}, we compute using again equation~\eqref{actsymmetry}
\begin{align*}
\mathrm{d}(\gamma Z)&=A\,\mathrm{d}Z(CZ+D)^{-1}-(AZ+B)(CZ+D)^{-1}C\,\mathrm{d}Z(CZ+D)^{-1} \\
&=(CZ+D)^{-t}((CZ+D)^{t}A\,\mathrm{d}Z-(AZ+B)^{t}C\,\mathrm{d}Z)(CZ+D)^{-1}.
\end{align*}
Using once again equation~\eqref{symprelation 1}, this gives
\begin{align}
\label{difZ calc2}
\mathrm{d}(\gamma Z)=(CZ+D)^{-t}\,\mathrm{d}Z\,(CZ+D)^{-1}.
\end{align}

The arclength $\mathrm{d}s^{2}_{n}$ and the volume form $\mathrm{d}\mu_{n}$ on $\mathbb{H}_{n}$ in terms of $Z=(z_{j,k})_{1\leq j\leq k\leq n}\in 
\mathbb{H}_{n}$ are given by
\begin{align*}
\mathrm{d}s^{2}_{n}(Z)&=\tr(Y^{-1}\,\mathrm{d}Z\,Y^{-1}\,\mathrm{d}\overline{Z})\quad\qquad(Z=X+iY), \\[1mm]
\mathrm{d}\mu_{n}(Z)&=\frac{\bigwedge\limits_{1\leq j\leq k\leq n}\mathrm{d}x_{j,k}\wedge\mathrm{d}y_{j,k}}{\det(Y)^{n+1}}\qquad(z_{j,k}=x_{j,k}+
iy_{j,k}).
\end{align*}
From equations~\eqref{imaginary transform} and~\eqref{difZ calc2} it is obvious that the arclength $\mathrm{d}s^{2}_{n}$ and the volume form 
$\mathrm{d}\mu_{n}$ on $\mathbb{H}_{n}$ given by the above equations are invariant under the symplectic action. Corresponding to this metric, 
we have the Laplace--Beltrami operator $\Delta$ on $\mathbb{H}_{n}$, called the Siegel Laplacian, which is also invariant under the symplectic 
action.  

\begin{definition}
Let $\Gamma\subset\mathrm{Sp}_{n}(\mathbb{R})$ be a subgroup commensurable with $\mathrm{Sp}_{n}(\mathbb{Z})$, i.e., the intersection
$\Gamma\cap\mathrm{Sp}_{n}(\mathbb{Z})$ is a finite index subgroup of $\Gamma$ as well as of $\mathrm{Sp}_{n}(\mathbb{Z})$. We let $\gamma_
{j}\in\mathrm{Sp}_{n}(\mathbb{Z})$ ($j=1,\ldots,h$) denote a set of representatives for the left cosets of $\Gamma\cap\mathrm{Sp}_{n}(\mathbb{Z})$
in $\mathrm{Sp}_{n}(\mathbb{Z})$. Then, a \emph{Siegel modular form of weight $k$ and degree $n$ for $\Gamma$} is a function $f\colon\mathbb
{H}_{n}\longrightarrow\mathbb{C}$ satisfying the following conditions:
\begin{itemize}
\item[(i)] 
$f$ is holomorphic;
\item[(ii)] 
$f(\gamma Z)=\det(CZ+D)^{k}f(Z)$ for all $\gamma=\big(\begin{smallmatrix}A&B\\C&D\end{smallmatrix}\big)\in\Gamma$;
\item[(iii)] 
given $Y_{0}\in\mathrm{Sym}_{n}(\mathbb{R})$ with $Y_{0}>0$, the quantities $\det(C_{j}Z+D_{j})^{-k}f(\gamma_{j} Z)$ are bounded in the region 
$\{Z=X+iY\in\mathbb{H}_{n}\,\vert\,Y\geq Y_{0}\}$ for the set of representatives $\gamma_{j}=\big(\begin{smallmatrix}A_{j}&B_{j}\\C_{j}&D_{j}\end
{smallmatrix}\big)\in\mathrm{Sp}_{n}(\mathbb{Z})$ ($j=1,\ldots,h$).
\end{itemize}
Moreover, a Siegel modular form $f$ as above is called a \emph{Siegel cusp form of weight $k$ and degree $n$ for $\Gamma$} if condition~(iii)
above is strengthened to the condition
\begin{itemize}
\item[(iii')] 
given $Y_{0}\in\mathrm{Sym}_{n}(\mathbb{R})$ with $Y_{0}\gg 0$, the quantities $\det(C_{j}Z+D_{j})^{-k}f(\gamma_{j} Z)$ become arbitrarily small
in the region $\{Z=X+iY\in\mathbb{H}_{n}\,\vert\,Y\geq Y_{0}\}$ for the set of representatives $\gamma_{j}=\big(\begin{smallmatrix}A_{j}&B_{j}\\C_
{j}&D_{j}\end{smallmatrix}\big)\in\mathrm{Sp}_{n}(\mathbb{Z})$ ($j=1,\ldots,h$).
\end{itemize}
\end{definition}

\begin{remark}
The sets of Siegel modular forms and Siegel cusp forms of weight $k$ and degree $n$ for $\Gamma$ have the structure of $\mathbb{C}$-vector 
spaces, which we denote by $\mathcal{M}_{k}^{n}(\Gamma)$ and $\mathcal{S}_{k}^{n}(\Gamma)$, respectively, and which turn out to be finite 
dimensional. Moreover, the space $\mathcal{S}_{k}^{n}(\Gamma)$ is equipped with the so-called Petersson inner product given by
\begin{align*}
\langle f,g\rangle:=\int_{\Gamma\backslash\mathbb{H}_{n}}\det(Y)^{k}f(Z)\overline{g}(Z)\,\mathrm{d}\mu_{n}(Z)\qquad(f,g\in\mathcal{S}_{k}^{n}
(\Gamma)),
\end{align*}
making $\mathcal{S}_{k}^{n}(\Gamma)$  into a hermitian inner product space.
\end{remark}

\section{Siegel--Maa\ss{} Laplacian of weight $(\alpha,\beta)$}
In this section, we will recall from~\cite{Maass2} various differential operators acting on smooth complex valued functions defined on $\mathbb{H}_
{n}$. In particular, we will define the Siegel--Maa\ss{} Laplacian of weight $(\alpha,\beta)$, where $\alpha,\beta\in\mathbb{R}$. Letting $\alpha=k/2$
and $\beta=-k/2$ will then lead us to the Siegel--Maa\ss{} Laplacian $\Delta_{k}$ mentioned in formula \eqref{deltak} in the introduction. We point
out that the Siegel Laplacian $\Delta$ mentioned in the previous section and the Siegel--Maa\ss{} Laplacian $\Delta_{k}$ are related by the formula
\begin{align*}
\Delta_{k}=\Delta-ik\tr\bigg(Y\dpd{}{X}\bigg)
\end{align*}
with the symmetric $(n\times n)$-matrix $\partial/\partial X$ of partial derivatives being defined below.

Given $Z=X+iY\in\mathbb{H}_{n}$, we start by introducing the following symmetric $(n\times n)$-matrices of partial derivatives:
\begin{align*}
\mathrm{(i)}&\quad\bigg(\dpd{}{X}\bigg)_{j,k}:=\frac{1+\delta_{j,k}}{2}\dpd{}{x_{j,k}}, \\
\mathrm{(ii)}&\quad\bigg(\dpd{}{Y}\bigg)_{j,k}:=\frac{1+\delta_{j,k}}{2}\dpd{}{y_{j,k}}, \\
\mathrm{(iii)}&\quad\dpd{}{Z}:=\frac{1}{2}\bigg(\dpd{}{X}-i\dpd{}{Y}\bigg), \\
\mathrm{(iv)}&\quad\dpd{}{\overline{Z}}:=\frac{1}{2}\bigg(\dpd{}{X}+i\dpd{}{Y}\bigg),
\end{align*}
where $\delta_{j,k}$ is the Kornecker delta symbol.

\begin{definition}
\label{maass operators}
Following Maa\ss{}~\cite{Maass2}, we define, using the above notations, for arbitrary real numbers $\alpha,\beta\in\mathbb{R}$, the following $(n
\times n)$-matrices of differential operators acting on smooth complex valued functions on~$\mathbb{H}_{n}$:
\begin{align*}
\mathrm{(i)}&\quad K_{\alpha}:=(Z-\overline{Z})\dpd{}{Z}+\alpha\mathbbm{1}_{n}, \\
\mathrm{(ii)}&\quad\Lambda_{\beta}:=(Z-\overline{Z})\dpd{}{\overline{Z}}-\beta\mathbbm{1}_{n}, \\[2mm]
\mathrm{(iii)}&\quad\Omega_{\alpha,\beta}:=\Lambda_{\beta-(n+1)/2}K_{\alpha}+\alpha(\beta-(n+1)/2)\mathbbm{1}_{n}, \\[3mm]
\mathrm{(iv)}&\quad\widetilde{\Omega}_{\alpha,\beta}:=K_{\alpha-(n+1)/2}\Lambda_{\beta}+\beta(\alpha-(n+1)/2)\mathbbm{1}_{n}.
\end{align*}
\end{definition}
Next, we want to expand $\Omega_{\alpha,\beta}$ and $\widetilde{\Omega}_{\alpha,\beta}$ in terms of $Z,\overline{ Z},\partial/\partial Z$, and 
$\partial/\partial{\overline{Z}}$. For that we need the following lemma.

\begin{lemma}
\label{operator derivative identity general}
Let $C,D\colon\mathbb{H}_{n}\longrightarrow\mathrm{M}_{n}(\mathbb{C})$ be smooth matrix valued functions depending on $Z$ and $\overline{Z}$. 
Then, the following matrix operator identities hold: \\[1mm]
\emph{(i)} Assuming that $\partial C/\partial Z=0$ and $\partial D/\partial Z=0$, we have
\begin{align*}
\dpd{}{Z}(CZ+D)^{t}=\bigg((CZ+D)\dpd{}{Z}\bigg)^{t}+\frac{1}{2}(n+1)C^{t}.
\end{align*}
\emph{(ii)} Assuming that $\partial C/\partial\overline{Z}=0$ and $\partial D/\partial\overline{Z}=0$, we have
\begin{align*}
\dpd{}{\overline{Z}}(C\overline{Z}+D)^{t}=\bigg((C\overline{Z}+D)\dpd{}{\overline{Z}}\bigg)^{t}+\frac{1}{2}(n+1)C^{t}.
\end{align*}
\end{lemma}
\begin{proof}
Since part~(ii) follows from part~(i) by conjugation, we prove only (i). Let $\Phi$ be a matrix depending on $Z$ and $\overline{Z}$ such that the
product $(CZ+D)^{t}\Phi$ makes sense. Then, writing the $(j,k)$-th entry of the matrix $\partial/\partial Z(CZ+D)^{t}\Phi$ as the sum
\begin{align*}
\bigg(\dpd{}{Z}(CZ+D)^{t}\Phi\bigg)_{j,k}=\sum\limits_{l,m=1}^{n}\bigg(\dpd{}{Z}\bigg)_{j,l}\big((CZ+D)^{t}_{l,m}\Phi_{m,k}\big)
\end{align*}
and noting that $\partial Z/\partial z_{j,l}=(1-\delta_{j,l})E_{j,l}+E_{l,j}$, where $E_{j,k}\in\mathrm{M}_{n}(\mathbb{C})$ is the matrix with its $(j,k)$-th 
entry being $1$ and the remaining entries being $0$, elementary calculations lead us to the operator identity
\begin{align*}
\dpd{}{Z}(CZ+D)^{t}=\bigg((CZ+D)\dpd{}{Z}\bigg)^{t}+\frac{1}{2}(n+1)C^{t},
\end{align*}
which is what we needed to prove.
\end{proof}

\begin{corollary}
\label{operator derivative identity}
For $Z\in\mathbb{H}_{n}$, the following operator identities hold:
\begin{align*}
\mathrm{(i)}\quad\dpd{}{Z}(Z-\overline{Z})=\bigg((Z-\overline{Z})\dpd{}{Z}\bigg)^{t}+\frac{1}{2}(n+1)\mathbbm{1}_{n}, \\
\mathrm{(ii)}\quad\dpd{}{\overline{Z}}(Z-\overline{Z})=\bigg((Z-\overline{Z})\dpd{}{\overline{Z}}\bigg)^{t}-\frac{1}{2}(n+1)\mathbbm{1}_{n}.
\end{align*}
\end{corollary}
\begin{proof}
As $\partial\overline{Z}/\partial Z=0$, we can choose $C=\mathbbm{1}_{n}$ and $D=-\overline{ Z}$ in Lemma~\ref{operator derivative identity 
general}~(i), from which the first claimed formula follows. The second formula follows analogously.
\end{proof}

Using the above corollary, one can expand $\Omega_{\alpha,\beta}$ and $\widetilde{\Omega}_{\alpha,\beta}$ as
\begin{align*}
\Omega_{\alpha,\beta}&=(Z-\overline{ Z})\bigg((Z-\overline{ Z})\dpd{}{\overline{Z}}\bigg)^{t}\dpd{}{Z}+\alpha(Z-\overline{Z})\dpd{}{\overline{Z}}-
\beta(Z-\overline{Z})\dpd{}{Z}, \\
\widetilde{\Omega}_{\alpha,\beta}&=(Z-\overline{Z})\bigg((Z-\overline{Z})\dpd{}{Z}\bigg)^{t}\dpd{}{\overline{Z}}+\alpha(Z-\overline{Z})\dpd{}
{\overline{Z}}-\beta(Z-\overline{Z})\dpd{}{Z}.
\end{align*}
Then, $\Omega_{\alpha,\beta}$ and $\widetilde{\Omega}_{\alpha,\beta}$ are related by the identity
\begin{align*}
\widetilde{\Omega}_{\alpha,\beta}=(Z-\overline{ Z})\big((Z-\overline{Z})^{-1}\Omega_{\alpha,\beta}\big)^{t}. 
\end{align*}
\begin{definition}
The operator $\Delta_{\alpha,\beta}:=-\tr(\Omega_{\alpha,\beta})=-\tr(\widetilde{\Omega}_{\alpha,\beta})$ is called the \emph{Siegel--Maa\ss{} 
Laplacian of weight $(\alpha,\beta)$}.
\end{definition}

\section{Transformation behaviour of Maa\ss{} operators}
Recall that the symplectic action of $\gamma=\big(\begin{smallmatrix}A&B\\C&D\end{smallmatrix}\big)\in\mathrm{Sp}_{n}(\mathbb{R})$ on the 
point $Z\in\mathbb{H}_{n}$ is given by
\begin{align*}
\gamma Z=(AZ+B)(CZ+D)^{-1}=(CZ+D)^{-t}(AZ+B)^{t};
\end{align*}
to avoid cumbersome notation, we will use in this paper sometimes the shorthand $Z^{\gamma}:=\gamma Z$. In this section, we will then study 
the transformation behaviour of the Maa\ss{} operators introduced in Definition~\ref{maass operators} by expressing the operators $K^{\gamma}_
{\alpha}$, $\Lambda^{\gamma}_{\beta}$, $\Omega^{\gamma}_{\alpha,\beta}$ obtained by replacing $Z,\overline{Z}$ by $Z^{\gamma},\overline
{Z}^{\gamma}$ in ${K}_{\alpha}$, ${\Lambda}_{\beta}$, ${\Omega}_{\alpha,\beta}$, respectively, as they operate on smooth complex valued 
functions defined on $\mathbb{H}_{n}$. 

We begin by investigating how the matrix differential operators $\partial/\partial Z$ and $\partial/\partial\overline{Z}$ transform under the symplectic 
action of $\gamma$ on $Z$. From equation~\eqref{difZ calc2}, we know that the differential $\mathrm{d}Z$ transforms like
\begin{align*}
\mathrm{d}Z^{\gamma}=(CZ+D)^{-t}\,\mathrm{d}Z\,(CZ+D)^{-1}. 
\end{align*}
Therefore, as the differential of a smooth function $\varphi\colon\mathbb{H}_{n}\longrightarrow\mathbb{C}$ depending only on $Z$ is given by 
$\mathrm{d}\varphi=\tr(\partial\varphi/\partial Z\,\mathrm{d}Z)$, we have
\begin{align*}
\dpd{\varphi}{Z}=(CZ+D)^{-1}\dpd{\varphi}{{Z^{\gamma}}}(C Z+D)^{-t},
\end{align*}
i.e., the operator $\partial/\partial Z$ transforms as
\begin{align}
\label{partial transform}
\dpd{}{{Z^{\gamma}}}=(C Z+D)\bigg((CZ+D)\dpd{}{Z}\bigg)^{t}.
\end{align}
By conjugation, the operator $\partial/\partial\overline{Z}$ transforms as
\begin{align}
\label{partial transform 2}
\dpd{}{{\overline{Z}^{\gamma}}}=(C\overline{Z}+D)\bigg((C\overline{Z}+D)\dpd{}{\overline{Z}}\bigg)^{t}.
\end{align}

Next we need to know how to differentiate $\det(Z-\overline{Z})$ and $\det(CZ+D)$ with respect to $Z$, which we carry out in the next two lemmas.
\begin{lemma}
\label{determinant differentiate}
The matrix identity 
\begin{align*}
\dpd{\det(Z-\overline{Z})}{Z}=\det(Z-\overline{Z})(Z-\overline{Z})^{-1}
\end{align*}
holds.
\end{lemma}
\begin{proof}
Since $Y=\frac{1}{2i}(Z-\overline{Z})\in\mathrm{Sym}_{n}(\mathbb{R})$, it can be diagonalized with orthogonal matrices. Thus, let $Y=U^{t}\Lambda 
U$, where $\Lambda\in\mathrm{M}_{n}(\mathbb{R})$ is a diagonal matrix with diagonal entries equal to the eigenvalues $\lambda_{1},\ldots,\lambda_
{n}$ of $Y$ and $U\in\mathrm{O}_{n}(\mathbb{R})$. Therefore, differentiating $\det(Y)$ with respect to the entries of $Y=(y_{j,k})_{1\leq j,k\leq n}$, 
we have
\begin{align*}
\dpd{\det(Y)}{y_{j,k}}=\dpd{\prod\limits_{l=1}^{n}\lambda_{l}}{y_{j,k}}=\det(Y)\tr\bigg(\Lambda^{-1}\dpd{\Lambda}{y_{j,k}}\bigg).
\end{align*}
Now writing $\tr(\Lambda^{-1}\,\partial\Lambda/\partial y_{j,k})$ as $\tr(Y^{-1}\,\partial Y/\partial y_{j,k})$ and using the fact that
\begin{align*}
\dpd{Y}{y_{j,k}}=(1-\delta_{j,k})E_{j,k}+E_{k,j}, 
\end{align*}
we obtain
\begin{align*}
\dpd{\det(Y)}{y_{j,k}}=\det (Y)(2-\delta_{j,k})(Y^{-1})_{j,k}. 
\end{align*}
Now since $((1+\delta_{j,k})/2)(2-\delta_{j,k})=1$, we have $\partial\det(Y)/\partial Y=\det(Y)Y^{-1}$. This is equivalent to the identity 
\begin{align*}
\dpd{\det(Z-\overline{Z})}{Z}=\det(Z-\overline{Z})(Z-\overline{Z})^{-1},
\end{align*}
which is what we needed to prove. 
\end{proof}

\begin{lemma}
\label{gen determinant differentiate}
Let $\gamma=\big(\begin{smallmatrix}A&B\\C&D\end{smallmatrix}\big)\in\mathrm{Sp}_{n}(\mathbb{R})$. Then, the matrix identity
\begin{align*}
\dpd{\det(CZ+D)}{Z}=\det(CZ+D)(CZ+D)^{-1}C=\det(CZ+D)C^{t}(C Z+D)^{-t}
\end{align*}
holds.
\end{lemma}
\begin{proof}
For $1\leq j,k\leq n$, let $u_{j,k}\colon\mathbb{H}_{n}\longrightarrow\mathbb{C}$ be smooth scalar valued functions; then, $U:=(u_{j,k})_{1\leq j,
k\leq n}$ becomes a smooth matrix valued function on $\mathbb{H}_{n}$. Moreover, let $\varphi\colon\mathrm{M}_{n}(\mathbb{C})\longrightarrow
\mathbb{C}$ be a smooth scalar valued function. Then, differentiating the function $\varphi\circ U\colon\mathbb{H}_{n}\longrightarrow\mathbb{C}$ 
with respect to the entries $z_{j,k}$ of $Z$, we have
\begin{align*}
\dpd{\varphi(U(Z))}{z_{j,k}}=\sum\limits_{l,m=1}^{n}\dpd{\varphi}{u_{l,m}}\dpd{u_{l,m}}{z_{j,k}}=\sum\limits_{l,m=1}^{n}\bigg(\dpd{\varphi}{U}\bigg)^
{t}_{m,l}\bigg(\dpd{U}{z_{j,k}}\bigg)_{l,m}=\sum\limits_{m=1}^{n}\bigg(\bigg(\dpd{\varphi}{U}\bigg)^{t}\dpd{U}{z_{j,k}}\bigg)_{m,m}.
\end{align*}
Thus, the chain rule of differentiation in this case takes the form
\begin{align*}
\dpd{\varphi(U(Z))}{ z_{j,k}}=\tr\bigg(\bigg(\dpd{\varphi}{U}\bigg)^t\dpd{U}{ z_{j,k}}\bigg).
\end{align*}
Note that since we do not assume $U$ to be symmetric beforehand, in this case we have $(\partial/\partial U)_{j,k}=\partial/\partial u_{j,k}$ instead 
of $(\partial/\partial U)_{j,k}=((1+\delta_{j,k})/2)\partial/\partial u_{j,k}$.

Now putting $\varphi(Z)=\det(Z)$ and $U(Z)=CZ+D$, by the above formula, we have
\begin{align*}
\dpd{\det(CZ+D)}{z_{j,k}}=\tr\bigg(\bigg(\dpd{\det(CZ+D)}{(CZ+D)}\bigg)^{t}\dpd{(CZ+D)}{z_{j,k}}\bigg). 
\end{align*}
We already derived a formula for differentiating the determinant of a symmetric matrix by its entries. The structure of symmetry actually complicates 
the calculation as its entries are no longer independent. For a square matrix $U$, not assumed to be symmetric beforehand, the calculation can be 
simplified by considering the cofactor expansion of a determinant. Let $\widetilde{U}=(\widetilde{u}_{j,l})_{1\leq j,l\leq n}$ be the cofactor matrix of 
$U$. Then, we have
\begin{align*}
\dpd{\det(U)}{u_{j,k}}=\dpd{}{u_{j,k}}\sum\limits_{l=1}^{n}u_{j,l}\widetilde{u}_{j,l}=\widetilde{u}_{j,k}.
\end{align*}
Here we exploit the property that since we delete the $j$-th row (and the $l$-th column) to build the cofactor $\widetilde{u}_{j,l}$, it must be independent 
of $u_{j,k}$. This does not hold for a symmetric matrix. Now, since $U^{-1}=1/\det (U)\,\widetilde{U}^{t}$, we have
\begin{align*}
\dpd{\det(U)}{U}=\widetilde{U}=\det(U)\,U^{-t}.
\end{align*}
Thus, going back to our initial calculation, we arrive at
\begin{align*}
\dpd{\det(CZ+D)}{z_{j,k}}&=\det(CZ+D)\tr\bigg((CZ+D)^{-1}\dpd{(CZ+D)}{z_{j,k}}\bigg) \\
&=\det(CZ+D)\sum\limits_{l,m=1}^{n}(CZ+D)^{-1}_{l,m}\dpd{(CZ+D)_{m,l}}{z_{j,k}}. 
\end{align*}
Now, entrywise partial differentiation with respect to entries $z_{j,k}$ of $Z$ followed by an elementary calculation with taking care of the ensuing 
Kronecker delta symbols leads us to the matrix identity
\begin{align*}
\dpd{\det(CZ+D)}{Z}=\det(CZ+D)(CZ+D)^{-1}C=\det(CZ+D)C^{t}(CZ+D)^{-t},
\end{align*}
which is what we needed to prove. 
\end{proof}

Lemmas~\ref{determinant differentiate}  and~\ref{gen determinant differentiate} prepare the groundwork for calculating the transformation behaviour 
of the Maa\ss{} operators, which we undertake one by one in the subsequent three propositions.

\begin{proposition}
\label{Ktransformation}
Let $\gamma=\big(\begin{smallmatrix}A&B\\C&D\end{smallmatrix}\big)\in\mathrm{Sp}_{n}(\mathbb{R})$ and $\varphi\colon\mathbb{H}_{n}\longrightarrow 
\mathbb{C}$ be a smooth function. Then, the operator $K^{\gamma}_{\alpha}$ obtained by replacing $Z\in\mathbb{H}_{n}$ in ${K}_{\alpha}$ by $Z^
{\gamma}=\gamma Z$ is related to the operator ${K}_{\alpha}$ by the identity
\begin{align*}
&K^{\gamma}_{\alpha}\big(\det(CZ+D)^{\alpha}\det(C\overline{Z}+D)^{\beta}\varphi(Z)\big) \\[1mm]
&\qquad=\det(CZ+D)^{\alpha}\det(C\overline{Z}+D)^{\beta}(C\overline{Z}+D)^{-t}K_{\alpha}\varphi(Z)(CZ+D)^{t}.
\end{align*}
\end{proposition}
\begin{proof}
From the definition of $K^{\gamma}_{\alpha}$, we have
\begin{align*}
&K^{\gamma}_{\alpha}\big(\det(CZ+D)^{\alpha}\det(C\overline{Z}+D)^{\beta}\varphi(Z)\big) \\
&\qquad=\bigg((Z^{\gamma}-\overline{Z}^{\gamma})\frac{\partial}{\partial Z^{\gamma}}+\alpha\mathbbm{1}_{n}\bigg)\det(CZ+D)^{\alpha}\det(C\overline{Z}+
D)^{\beta}\varphi(Z).
\end{align*}
Then, expanding $\partial/\partial Z^{\gamma}$ by means of equation~\eqref{partial transform} gives
\begin{align}
\label{longKstar}
\notag
&K^{\gamma}_{\alpha}\big(\det(CZ+D)^{\alpha}\det(C\overline{Z}+D)^{\beta}\varphi(Z)\big)=\alpha\det(CZ+D)^{\alpha}\det(C\overline{Z}+D)^{\beta}\varphi
(Z)\mathbbm{1}_{n} \\[1mm]
&\qquad+(C\overline{Z}+D)^{-t}(Z-\overline{Z})\dpd{}{Z}\big(\det(CZ+D)^{\alpha}\det(C\overline{Z}+D)^{\beta}\varphi(Z)\big)(CZ+D)^{t}.
\end{align}
Now, focusing on the second line of the above equality and using Lemma~\ref{gen determinant differentiate}, we get
\begin{align*}
&\dpd{}{Z}\big(\det(CZ+D)^{\alpha}\det(C\overline{Z}+D)^{\beta}\varphi(Z)\big) \\[1mm]
&\qquad=\det(CZ+D)^{\alpha}\det(C\overline{Z}+D)^{\beta}\bigg(\alpha\varphi(Z)C^{t}(CZ+D)^{-t}+\dpd{\varphi}{Z}\bigg).
\end{align*}
Multiplying the above equation from the left by $(Z-\overline{Z})$ gives
\begin{align*}
&(Z-\overline{Z})\dpd{}{Z}\big(\det(CZ+D)^{\alpha}\det(C\overline{Z}+D)^{\beta}\varphi(Z)\big) \\[1mm]
&\qquad=\det(CZ+D)^{\alpha}\det(C\overline{Z}+D)^{\beta}\bigg(\alpha\varphi(Z)(ZC^{t}-\overline{Z}C^{t})(CZ+D)^{-t}+(Z-\overline{Z})\dpd{\varphi}{Z}\bigg).
\end{align*}
Now writing $(ZC^{t}-\overline{Z}C^{t})=(CZ+D)^{t}-(C\overline{Z}+D)^{t}$ and using the definition of $K_{\alpha}$ on the right-hand side of the above equation, 
we have 
\begin{align*}
&(Z-\overline{Z})\dpd{}{Z}\big(\det(CZ+D)^{\alpha}\det(C\overline{Z}+D)^{\beta}\varphi(Z)\big) \\[1mm]
&\qquad=\det(C Z+D)^{\alpha}\det(C\overline{ Z}+D)^{\beta} \big(K_{\alpha}\varphi( Z)-\alpha\varphi( Z)(C\overline{ Z}+D)^t(C Z+D)^{-t}\big).
\end{align*}
Therefore, multiplying on the left by $(C\overline{Z}+D)^{-t}$ and on the right by $(CZ+D)^{t}$, we obtain
\begin{align*}
&(C\overline{Z}+D)^{-t}(Z-\overline{Z})\dpd{}{Z}\big(\det(CZ+D)^{\alpha}\det(C\overline{Z}+D)^{\beta}\varphi(Z)\big)(CZ+D)^{t} \\[1mm]
&\qquad=\det(CZ+D)^{\alpha}\det(C\overline{Z}+D)^{\beta}\big((C\overline{Z}+D)^{-t}K_{\alpha}\varphi(Z)(CZ+D)^{t}-\alpha\varphi(Z)\mathbbm{1}_{n}\big).
\end{align*}
Combining the last equality with equation~\eqref{longKstar}, leads to the identity
\begin{align*}
&K^{\gamma}_{\alpha}\big(\det(CZ+D)^{\alpha}\det(C\overline{Z}+D)^{\beta}\varphi(Z)\big) \\[1mm]
&\qquad=\det(CZ+D)^{\alpha}\det(C\overline{Z}+D)^{\beta}(C\overline{Z}+D)^{-t}{K}_{\alpha}\varphi(Z)(CZ+D)^{t},
\end{align*}
which is what we had set out to prove.
\end{proof}

\begin{proposition}
\label{Ltransformation}
Let $\gamma=\big(\begin{smallmatrix}A&B\\C&D\end{smallmatrix}\big)\in\mathrm{Sp}_{n}(\mathbb{R})$ and $\varphi\colon\mathbb{H}_{n}\longrightarrow 
\mathbb{C}$ be a smooth function. Then, the operator $\Lambda^{\gamma}_{\beta}$ obtained by replacing $Z\in\mathbb{H}_{n}$ in ${\Lambda}_{\beta}$ by 
$Z^{\gamma}=\gamma Z$ is related to the operator ${\Lambda}_{\beta}$ by the identity
\begin{align*}
&\Lambda^{\gamma}_{\beta}\big(\det(CZ+D)^{\alpha}\det(C\overline{Z}+D)^{\beta}\varphi(Z)\big) \\[1mm]
&\qquad=\det(CZ+D)^{\alpha}\det(C\overline{Z}+D)^{\beta}(CZ+D)^{-t}{\Lambda}_{\beta}\varphi(Z)(C\overline{Z}+D)^{t}. 
\end{align*}
\end{proposition}
\begin{proof}
Since $\overline{K}_{\beta}=-\Lambda_{\beta}$, the required identity follows from Proposition~\ref{Ktransformation} by complex conjugation.
\end{proof}

\begin{proposition}
\label{Omegatransform}
Let $\gamma=\big(\begin{smallmatrix}A&B\\C&D\end{smallmatrix}\big)\in\mathrm{Sp}_{n}(\mathbb{R})$ and $\varphi\colon\mathbb{H}_{n}\longrightarrow 
\mathbb{C}$ be a smooth function. Then, the operator $\Omega^{\gamma}_{\alpha,\beta}$ obtained by replacing $Z\in\mathbb{H}_{n}$ in ${\Omega}_{\alpha,
\beta}$ by $Z^{\gamma}=\gamma Z$ is related to the operator ${\Omega}_{\alpha,\beta}$ by the identity
\begin{align}
\label{Omegatransformeq}
\notag
&\Omega^{\gamma}_{\alpha,\beta}\big(\det(CZ+D)^{\alpha}\det(C\overline{Z}+D)^{\beta}\varphi(Z)\big) \\[1mm]
&\qquad=\det(CZ+D)^{\alpha}\det(C\overline{Z}+D)^{\beta}(CZ+D)^{-t}{\Omega}_{\alpha,\beta}\varphi(Z)(CZ+D)^{t}.
\end{align}
\end{proposition}
\begin{proof}
To prove the proposition, we first need to calculate 
\begin{align*}
\dpd{}{{\overline{Z}^{\gamma}}}\big(K^{\gamma}_{\alpha}\big(\det(CZ+D)^{\alpha}\det(C\overline{Z}+D)^{\beta}\varphi(Z)\big)\big),
\end{align*}
which, upon expanding $\partial/\partial\overline{Z}^{\gamma}$ by means of equation~\eqref{partial transform 2} and using Proposition~\ref{Ktransformation},
becomes
\begin{align*}
(C\overline{Z}+D)\bigg((C\overline{Z}+D)\dpd{}{\overline{Z}}\bigg)^{t}\big(\det(CZ+D)^{\alpha}\det(C\overline{Z}+D)^{\beta}(C\overline{Z}+D)^{-t}{K}_{\alpha}
\varphi(Z)(CZ+D)^{t}\big). 
\end{align*}
Using Lemma~\ref{gen determinant differentiate}, the above expression becomes
\begin{align*}
&\det(CZ+D)^{\alpha}\det(C\overline{Z}+D)^{\beta}\bigg(\beta(C\overline{Z}+D)C^{t}(C\overline{Z}+D)^{-t}K_{\alpha}\varphi(Z)(CZ+D)^{t} \\
&\qquad+(C\overline{Z}+D)\bigg((C\overline{Z}+D)\dpd{}{\overline{Z}}\bigg)^{t}\big((C\overline{Z}+D)^{-t}K_{\alpha}\varphi(Z)(CZ+D)^{t}\big)\bigg).
\end{align*}
Now, using Lemma~\ref{operator derivative identity general}~(ii), we have
\begin{align*}
& \bigg((C\overline{Z}+D)\dpd{}{\overline{ Z}}\bigg)^{t}\big((C\overline{Z}+D)^{-t}{K}_{\alpha}\varphi(Z)(CZ+D)^{t}\big) \\
&\qquad=\dpd{}{\overline{ Z}}\big({K}_{\alpha}\varphi(Z)\big)(CZ+D)^{t}-\frac{1}{2}(n+1)C^{t}(C\overline{Z}+D)^{-t}{K}_{\alpha}\varphi(Z)(CZ+D)^{t},
\end{align*}
and thus arrive from the above calculation at
\begin{align*}
&\dpd{}{{\overline{Z}^{\gamma}}}\big(K^{\gamma}_{\alpha}\big(\det(CZ+D)^{\alpha}\det(C\overline{Z}+D)^{\beta}\varphi(Z)\big)\big) \\
&\qquad=\det(CZ+D)^{\alpha}\det(C\overline{Z}+D)^{\beta}\bigg((\beta-(n+1)/2)(C\overline{Z}+D)C^{t}(C\overline{Z}+D)^{-t}\times \\ 
&\qquad\quad\times K_{\alpha}\varphi(Z)(CZ+D)^{t}+(C\overline{Z}+D)\dpd{}{\overline{Z}}\big(K_{\alpha}\varphi(Z)\big)(CZ+D)^{t}\bigg).
\end{align*}
Then, multiplying both sides from the left by $(Z^{\gamma}-\overline{Z}^{\gamma})$ and using 
\begin{align*}
(Z^{\gamma}-\overline{Z}^{\gamma})=(C{Z}+D)^{-t}(Z-\overline{Z})(C\overline{Z}+D)^{-1}
\end{align*}
on the right-hand side, we get 
\begin{align*}
&(Z^{\gamma}-\overline{Z}^{\gamma})\dpd{}{{\overline{Z}^{\gamma}}}\big(K^{\gamma}_{\alpha}\big(\det(CZ+D)^{\alpha}\det(C\overline{Z}+D)^{\beta}\varphi
(Z)\big)\big) \\
&\qquad=\det(CZ+D)^{\alpha}\det(C\overline{Z}+D)^{\beta}\bigg((\beta-(n+1)/2)(CZ+D)^{-t}(Z-\overline{Z})C^{t}(C\overline{Z}+D)^{-t}\times \\
&\qquad\quad\times K_{\alpha}\varphi(Z)(CZ+D)^{t}+(CZ+D)^{-t}(Z-\overline{Z})\dpd{}{\overline{Z}}\big({K}_{\alpha}\varphi(Z)\big)(CZ+D)^{t}\bigg).
\end{align*} 
Next, writing the expression $(Z-\overline{Z})C^{t}$ on the right-hand side of the above equation as $(CZ+D)^{t}-(C\overline{Z}+D)^{t}$, we can rewrite the 
above equation as
\begin{align*}
&(Z^{\gamma}-\overline{Z}^{\gamma})\dpd{}{{\overline{Z}^{\gamma}}}\big(K^{\gamma}_{\alpha}\big(\det(CZ+D)^{\alpha}\det(C\overline{Z}+D)^{\beta}\varphi
(Z)\big)\big)=\det(CZ+D)^{\alpha}\det(C\overline{Z}+D)^{\beta}\times \\
&\qquad\times\bigg((\beta-(n+1)/2)\big((C\overline{Z}+D)^{-t}K_{\alpha}\varphi(Z)(CZ+D)^{t}-(CZ+D)^{-t}K_{\alpha}\varphi(Z)(CZ+D)^{t}\big) \\
&\qquad\quad+(CZ+D)^{-t}(Z-\overline{Z})\dpd{}{\overline{Z}}\big(K_{\alpha}\varphi(Z)\big)(CZ+D)^{t}\bigg).
\end{align*}
Now shifting the first term on the right-hand side to the left and using the transformation behaviour of $K_{\alpha}$ derived in Proposition~\ref{Ktransformation}, 
we arrive at
\begin{align*}
&\bigg((Z^{\gamma}-\overline{Z}^{\gamma})\dpd{}{{\overline{Z}^{\gamma}}}-(\beta-(n+1)/2)\mathbbm{1}_{n}\bigg)\big(K^{\gamma}_{\alpha}\big(\det(CZ+D)^
{\alpha}\det(C\overline{Z}+D)^{\beta}\varphi(Z)\big)\big) \\
&\qquad=\det(CZ+D)^{\alpha}\det(C\overline{Z}+D)^{\beta}(CZ+D)^{-t} \bigg((Z-\overline{Z})\dpd{}{\overline{Z}}-(\beta-(n+1)/2)\mathbbm{1}_{n}\bigg)\times 
\\[1mm]
&\qquad\quad\times K_{\alpha}\varphi(Z)(CZ+D)^{t},
\end{align*}
which, by Definition~\ref{maass operators} and the transformation behaviour of $\Lambda_{\beta}$ given in Proposition~\ref{Ltransformation}, gives 
\begin{align*}
&\Lambda^{\gamma}_{\beta-(n+1)/2}K^{\gamma}_{\alpha}\big(\det(CZ+D)^{\alpha}\det(C\overline{Z}+D)^{\beta}\varphi(Z)\big) \\[1mm]
&\qquad=\det(CZ+D)^{\alpha}\det(C\overline{Z}+D)^{\beta}(CZ+D)^{-t}\Lambda_{\beta-(n+1)/2}K_{\alpha}\varphi(Z)(CZ+D)^{t},
\end{align*}
which, by definition of $\Omega^{\gamma}_{\alpha,\beta}$, yields the desired identity.  
\end{proof}

\begin{remark}
\label{functional equation}
Let $\gamma=\big(\begin{smallmatrix}A&B\\C&D\end{smallmatrix}\big)\in\mathrm{Sp}_{n}(\mathbb{R})$ and $\varphi\colon\mathbb{H}_{n}\longrightarrow 
\mathbb{C}$ be a smooth function. Taking traces on both sides of equation~\eqref{Omegatransformeq} leads to the following transformation behaviour of 
the Siegel--Maa\ss{} Laplacian $\Delta_{\alpha,\beta}$
\begin{align*}
\Delta^{\gamma}_{\alpha,\beta}\big(\det(CZ+D)^{\alpha}\det(C\overline{Z}+D)^{\beta}\varphi(Z)\big)=\det(CZ+D)^{\alpha}\det(C\overline{Z}+D)^{\beta}\Delta_
{\alpha,\beta}\varphi(Z). 
\end{align*}
Now, if the smooth function $\varphi$ satisfies the functional equation
\begin{align*}
\varphi(Z^{\gamma})=\det(CZ+D)^{\alpha}\det(C\overline{Z}+D)^{\beta}\varphi(Z),
\end{align*}
the transformation behaviour of $\Delta_{\alpha,\beta}$ leads to the identity
\begin{align*}
\Delta^{\gamma}_{\alpha,\beta}\varphi(Z^{\gamma})=\det(CZ+D)^{\alpha}\det(C\overline{Z}+D)^{\beta}\Delta_{\alpha,\beta}\varphi(Z).
\end{align*}
\end{remark}

\begin{definition}
\label{siegel-maass definition}
Let $\Gamma\subset\mathrm{Sp}_{n}(\mathbb{R})$ be a subgroup commensurable with $\mathrm{Sp}_{n}(\mathbb{Z})$, i.e., the intersection $\Gamma\cap
\mathrm{Sp}_{n}(\mathbb{Z})$ is a finite index subgroup of $\Gamma$ as well as of $\mathrm{Sp}_{n}(\mathbb{Z})$. We let $\gamma_{j}\in\mathrm{Sp}_{n}
(\mathbb{Z})$ ($j=1,\ldots,h$) denote a set of representatives for the left cosets of $\Gamma\cap\mathrm{Sp}_{n}(\mathbb{Z})$ in $\mathrm{Sp}_{n}(\mathbb
{Z})$. We then let $\mathcal{V}_{\alpha,\beta}^{n}(\Gamma)$ denote the space of all functions $\varphi\colon\mathbb{H}_{n}\longrightarrow\mathbb{C}$ 
satisfying the following conditions:
\begin{itemize}
\item[(i)]
$\varphi$ is real-analytic;
\item[(ii)]
$\varphi(\gamma Z)=\det(CZ+D)^{\alpha}\det(C\overline{Z}+D)^{\beta}\varphi(Z)$ for all $\gamma=\big(\begin{smallmatrix}A&B\\C&D\end{smallmatrix}\big)\in
\Gamma$;
\item[(iii)]
given $Y_{0}\in\mathrm{Sym}_{n}(\mathbb{R})$ with $Y_{0}>0$, there exist $M\in\mathbb{R}_{>0}$ and $N\in\mathbb{N}$ such that the inequalities
\begin{align*}
\vert\det(C_{j}Z+D_{j})^{-\alpha}\det(C_{j}\overline{Z}+D_{j})^{-\beta}\varphi(\gamma_{j} Z)\vert\leq M\tr(Y)^{N} 
\end{align*}
hold in the region $\{Z=X+iY\in\mathbb{H}_{n}\,\vert\,Y\geq Y_{0}\}$ for the set of representatives $\gamma_{j}=\big(\begin{smallmatrix}A_{j}&B_{j}\\C_{j}&
D_{j}\end{smallmatrix}\big)\in\mathrm{Sp}_{n}(\mathbb{Z})$ ($j=1,\ldots,h$).
\end{itemize}
\end{definition}

\begin{remark}
\label{boundedness}
For $\varphi\in\mathcal{V}_{\alpha,\beta}^{n}(\Gamma)$, we set
\begin{align*}
\Vert{\varphi}\Vert^{2}:=\int\limits_{\Gamma\backslash\mathbb{H}_{n}}\det(Y)^{\alpha+\beta}\vert\varphi(Z)\vert^{2}\,\mathrm{d}\mu_{n}(Z),
\end{align*}
whenever it is defined. In this way we obtain the Hilbert space
\begin{align*}
\mathcal{H}_{\alpha,\beta}^{n}(\Gamma):=\big\{\varphi\in\mathcal{V}_{\alpha,\beta}^{n}(\Gamma)\,\big\vert\,\Vert\varphi\Vert<\infty\big\}
\end{align*}
equipped with the inner product
\begin{align*}
\langle\varphi,\psi\rangle=\int\limits_{\Gamma\backslash\mathbb{H}_{n}}\det(Y)^{\alpha+\beta}\varphi(Z)\overline{\psi}(Z)\,\mathrm{d}\mu_{n}(Z)\qquad(\varphi,
\psi\in\mathcal{H}_{\alpha,\beta}^{n}(\Gamma)).
\end{align*}
We note that in order to enable $\Vert\varphi\Vert<\infty$, the exponent $N\in\mathbb{N}$ in part (iii) of Definition~\ref{siegel-maass definition} has to be~$0$.
Moreover, we note that Remark~\ref{functional equation} shows that the Siegel--Maa\ss{} Laplacian $\Delta_{\alpha,\beta}$ acts on the Hilbert space $\mathcal
{H}_{\alpha,\beta}^{n}(\Gamma)$.
\end{remark}

\begin{definition}
Let $\Gamma\subset\mathrm{Sp}_{n}(\mathbb{R})$ be a subgroup commensurable with $\mathrm{Sp}_{n}(\mathbb{Z})$. The elements of the Hilbert space 
$\mathcal{H}_{\alpha,\beta}^{n}(\Gamma)$ are called \emph{automorphic forms of weight $(\alpha,\beta)$ and degree $n$ for $\Gamma$}. Moreover, if 
$\varphi\in\mathcal{H}_{\alpha,\beta}^{n}(\Gamma)$ is an eigenform of $\Delta_{\alpha,\beta}$, it is called a \emph{Siegel--Maa\ss{} form of weight $(\alpha,
\beta)$ and degree $n$ for $\Gamma$}.
\end{definition}

\begin{corollary}
\label{maass form transform}
Let $\Gamma\subset\mathrm{Sp}_{n}(\mathbb{R})$ be a subgroup commensurable with $\mathrm{Sp}_{n}(\mathbb{Z})$ and $\varphi\in\mathcal{H}_{\alpha,
\beta}^{n}(\Gamma)$. Then, we have for all $\gamma=\big(\begin{smallmatrix}A&B\\C&D\end{smallmatrix}\big)\in\Gamma$
\begin{align*}
\mathrm{(i)}&\quad K^{\gamma}_{\alpha}\varphi(Z^{\gamma})=\det(CZ+D)^{\alpha}\det(C\overline{Z}+D)^{\beta}(C\overline{Z}+D)^{-t}K_{\alpha}\varphi(Z)(CZ+D)
^{t}, \\[1mm]
\mathrm{(ii)}&\quad\Lambda^{\gamma}_{\beta}\varphi(Z^{\gamma})=\det(CZ+D)^{\alpha}\det(C\overline{Z}+D)^{\beta}(CZ+D)^{-t}{\Lambda}_{\beta}\varphi(Z)(C
\overline{Z}+D)^{t}, \\[1mm]
\mathrm{(iii)}&\quad\Omega^{\gamma}_{\alpha,\beta}\varphi(Z^{\gamma})=\det(CZ+D)^{\alpha}\det(C\overline{Z}+D)^{\beta}(CZ+D)^{-t}{\Omega}_{\alpha,\beta}
\varphi( Z)(CZ+D)^{t}.
\end{align*}
\end{corollary}
\begin{proof}
The proof is an immediate consequence of Propositions \ref{Ktransformation}--\ref{Omegatransform} and the definition of the Hilbert space $\mathcal{H}_{\alpha,
\beta}^{n}(\Gamma)$.
\end{proof}

\section{Symmetry of the Siegel--Maa\ss{} Laplacian of weight $(\alpha,\beta)$}
Let $\mathrm{d}Z:=(\mathrm{d}z_{j,k})_{1\leq j,k\leq n}$ denote the $(n\times n)$-matrix consisting of differential forms of degree $1$ and let $[\mathrm{d}Z]:=
\bigwedge_{1\leq j\leq k\leq n}\mathrm{d}z_{j,k}$ denote the differential form of degree $n(n+1)/2$ at $Z\in\mathbb{H}_{n}$. We introduce an $(n\times n)$-matrix $
\{\mathrm{d}Z\}$ consisting of differential forms of degree $(n(n+1)/2-1)$, namely
\begin{align*}
\{\mathrm{d}Z\}_{j,k}:=\frac{1+\delta_{j,k}}{2}\,\varpi_{j,k},
\end{align*}
where $\varpi_{j,k}$ is defined by
\begin{align*}
\varpi_{j,k}:=\varepsilon_{j,k}\bigwedge\limits_{\substack{1\leq l\leq m \leq n\\(l,m)\neq(j,k)}}\mathrm{d}z_{l,m}\qquad(1\leq j\leq k\leq n)
\end{align*} 
in case $j\leq k$ and $\varpi_{j,k}=\varpi_{k,j}$ in case $j>k$ with the sign $\varepsilon_{j,k}=\pm1$ determined by  $\mathrm{d}z_{j,k}\wedge\varpi_{j,k}=[\mathrm
{d}Z]$. It is easy to see that
\begin{align*}
\mathrm{d}Z\wedge\{\mathrm{d}Z\}=\frac{1}{2}(n+1)[\mathrm{d}Z]\mathbbm{1}_{n}.
\end{align*}
Let now $\gamma=\big(\begin{smallmatrix}A&B\\C&D\end{smallmatrix}\big)\in\mathrm{Sp}_{n}(\mathbb{R})$. Since we have $\mathrm{d}Z^{\gamma}=(CZ+D)^{-t}
\,\mathrm{d}Z\,(CZ+D)^{-1}$ and $[\mathrm{d}Z^{\gamma}]=\linebreak\det(CZ+D)^{-(n+1)}[\mathrm{d}Z]$, we derive from the relation
\begin{align*}
\mathrm{d}Z^{\gamma}\wedge\{\mathrm{d}Z^{\gamma}\}=\frac{1}{2}(n+1)[\mathrm{d}Z^{\gamma}]\mathbbm{1}_{n}
\end{align*} 
that the matrix $\{\mathrm{d}Z\}$ has the transformation behaviour
\begin{align*}
\{\mathrm{d}Z^{\gamma}\}=\det(C Z+D)^{-(n+1)}(CZ+D)\{\mathrm{d}Z\}(CZ+D)^{t}.
\end{align*}
Next we shall use these differential forms to show that the Siegel--Maa\ss{} Laplacian $\Delta_{\alpha,\beta}$ acts as a symmetric operator on the Hilbert space 
$\mathcal{H}_{\alpha,\beta}^{n}(\Gamma)$.

\begin{theorem}
\label{main theorem}
Let $\Gamma\subset\mathrm{Sp}_{n}(\mathbb{R})$ be a subgroup commensurable with $\mathrm{Sp}_{n}(\mathbb{Z})$ and let $\varphi,\psi\in\mathcal{H}_{\alpha,
\beta}^{n}(\Gamma)$ be compactly supported. Then, we have the formula
\begin{align*}
\langle-\Delta_{\alpha,\beta}\varphi,\psi\rangle=\int\limits_{\Gamma\backslash\mathbb{H}_{n}}\det(Y)^{\alpha+\beta}\tr\big(\Lambda_{\beta}\varphi(Z)\overline{\Lambda}_
{\beta}\overline{\psi}(Z)\big)\,\mathrm{d}{\mu}_{n}(Z)+n\beta(\alpha-(n+1)/2)\langle\varphi,\psi\rangle. 
\end{align*}
In particular, this formula establishes the relation
\begin{align*}
\langle\Delta_{\alpha,\beta}\varphi,\psi\rangle=\langle\varphi,\Delta_{\alpha,\beta}\psi\rangle,
\end{align*}
which shows that the Siegel--Maa\ss{} Laplacian $\Delta_{\alpha,\beta}$ acts as a symmetric operator on the Hilbert space $\mathcal{H}_{\alpha,\beta}^{n}(\Gamma)$.
\end{theorem}
\begin{proof}
We start by considering the differential form
\begin{align*}
\omega(Z):=\det(Z-\overline{Z})^{\alpha+\beta-(n+1)}\overline{\psi}(Z)\tr\big(\Lambda_{\beta}\varphi(Z)(Z-\overline{Z})\{\mathrm{d}Z\}\big)\wedge[\mathrm{d}\overline
{Z}].
\end{align*} 
Let $\gamma=\big(\begin{smallmatrix}A&B\\C&D\end{smallmatrix}\big)\in\Gamma$. Then, the transformation formulas
\begin{align*}
\mathrm{(a})&\quad\det(Z^{\gamma}-\overline{Z}^{\gamma})^{\alpha+\beta-(n+1)} \\
&\qquad=\det(CZ+D)^{-(\alpha+\beta-(n+1))}\det(C\overline{Z}+D)^{-(\alpha+\beta-(n+1))}\det(Z-\overline{Z})^{\alpha+\beta-(n+1)}, \\[1mm]
\mathrm{(b})&\quad\overline{\psi}(Z^{\gamma})=\det(CZ+D)^{\beta}\det(C\overline{Z}+D)^{\alpha}\overline{\psi}(Z), \\[1mm]
\mathrm{(c})&\quad\tr\big(\Lambda^{\gamma}_{\beta}\varphi(Z^{\gamma})(Z^{\gamma}-\overline{Z}^{\gamma})\{\mathrm{d}Z^{\gamma}\}\big) \\[1mm]
&\qquad=\det(CZ+D)^{\alpha-(n+1)}\det(C\overline{Z}+D)^{\beta}\tr\big(\Lambda_{\beta}\varphi(Z)(Z-\overline{Z})\{\mathrm{d}Z\}\big), \\[1mm]
\mathrm{(d})&\quad[\mathrm{d}\overline{Z}^{\gamma}]=\det(C\overline{Z}+D)^{-(n+1)}[\mathrm{d}\overline{Z}]
\end{align*}
show that $\omega(Z^{\gamma})=\omega(Z)$ for all $\gamma\in\Gamma$, i.e., $\omega(Z)$ is a $\Gamma$-invariant differential form on $\mathbb{H}_{n}$, and 
hence can be considered as a differential form on the quotient space $\Gamma\backslash\mathbb{H}_{n}$. Since the automorphic forms $\varphi,\psi$ are
real-analytic, the differential form $\omega$ is a smooth differential form. Therefore, by Stokes' theorem, we have
\begin{align*}
\int\limits_{\Gamma\backslash\mathbb{H}_{n}}\mathrm{d}\omega(Z)=\int\limits_{\partial\Gamma\backslash\mathbb{H}_{n}}\omega(Z).
\end{align*}
As $\varphi,\psi$ are compactly supported, the integral on the right-hand side of the above equation vanishes, which gives
\begin{align}
\label{int zero}
\int\limits_{\Gamma\backslash\mathbb{H}_{n}}\mathrm{d}\omega(Z)=0. 
\end{align}
As we shall see, by explicitly computing $\mathrm{d}\omega(Z)$, the vanishing of the above integral will lead to the formula claimed in the theorem.

For the computation of $\mathrm{d}\omega(Z)$, we set $\rho:=\det(Z-\overline{Z})^{\alpha+\beta-(n+1)}\overline{\psi}(Z)$, $P:=\Lambda_{\beta}\varphi(Z)$, and 
$Q:=(Z-\overline{Z})$. Then, we obtain
\begin{align*}
\omega(Z)=\rho\tr(P\,Q\,\{\mathrm{d}Z\})\wedge[\mathrm{d}\overline{Z}]=\sum\limits_{j,k,l=1}^{n}\rho\,p_{j,k}\,q_{k,l}\,\{\mathrm{d}Z\}_{l,j}\wedge[\mathrm{d}\overline
{Z}].
\end{align*}
Taking exterior derivatives on both sides leads to
\begin{align}
\notag
\mathrm{d}\omega(Z)&=\sum\limits_{j,k,l=1}^{n}\dpd{}{z_{l,j}}(\rho\,p_{j,k}\,q_{k,l})\,\mathrm{d}z_{l,j}\wedge\frac{1+\delta_{l,j}}{2}\varpi_{l,j}\wedge[\mathrm{d}\overline
{Z}] \\
\notag
&=\sum\limits_{j,k,l=1}^{n}\frac{1+\delta_{l,j}}{2}\dpd{}{z_{l,j}}(\rho\,p_{j,k}\,q_{k,l})\,[\mathrm{d}Z]\wedge[\mathrm{d}\overline{Z}] \\
\label{dabba}
&=\sum\limits_{j,k,l=1}^{n}\bigg(\dpd{}{Z}\bigg)_{l,j}(\rho\,p_{j,k}\,q_{k,l})\,[\mathrm{d}Z]\wedge[\mathrm{d}\overline{Z}].
\end{align}
Now a term by term differentiation in the last expression on the right-hand side of the above equation allows us to write it as the sum of the three traces
\begin{align}
\label{rhopq0}
\sum\limits_{j,k,l=1}^{n}\bigg(\dpd{}{Z}\bigg)_{l,j}(\rho\,p_{j,k}\,q_{k,l})=\tr\bigg(\dpd{\rho}{Z}P\,Q\bigg)+\rho\tr\bigg(\dpd{}{Z}P\,Q\bigg)+\rho\tr\bigg(P^t\dpd{}{Z}Q\bigg),
\end{align}
which we calculate one by one next.

\noindent
(i) We begin by considering 
\begin{align*}
\dpd{\rho}{Z}=\dpd{}{Z}\big(\det(Z-\overline{Z})^{\alpha+\beta-(n+1)}\overline{\psi}(Z)\big),
\end{align*}
which, by Lemma~\ref{determinant differentiate}, calculates to
\begin{align*}
\dpd{\rho}{Z}=(\alpha+\beta-(n+1))\det(Z-\overline{Z})^{\alpha+\beta-(n+1)}(Z-\overline{Z})^{-1}\,\overline{\psi}(Z)+\det(Z-\overline{Z})^{\alpha+\beta-(n+1)}\dpd
{\overline{\psi}(Z)}{Z}.
\end{align*}
Now multiplying both sides of the above equation on the right by $P\,Q=\Lambda_{\beta}\varphi(Z)(Z-\overline{Z})$ and taking the trace gives
\begin{align*}
\tr\bigg(\dpd{\rho}{Z}P\,Q\bigg)&=\det(Z-\overline{Z})^{\alpha+\beta-(n+1)}\bigg((\alpha+\beta-(n+1))\tr\big((Z-\overline{Z})^{-1}\overline{\psi}(Z)\Lambda_{\beta}
\varphi(Z)(Z-\overline{Z})\big) \\
&\hspace*{35mm}+\tr\bigg(\dpd{\overline{\psi}(Z)}{Z}\Lambda_{\beta}\varphi(Z)(Z-\overline{Z})\bigg)\bigg),
\end{align*}
which, upon rearranging the terms inside the traces on the right-hand side by cyclically permuting them, becomes 
\begin{align}
\notag
&\tr\bigg(\dpd{\rho}{Z}P\,Q\bigg)=\det(Z-\overline{Z})^{\alpha+\beta-(n+1)}\bigg((\alpha+\beta-(n+1))\tr\big(\Lambda_{\beta}\varphi(Z)\overline{\psi}(Z)\big) \\ 
\label{rhopq1}
&\hspace*{57mm}+\tr\bigg(\Lambda_{\beta}\varphi(Z)(Z-\overline{Z})\dpd{\overline{\psi}(Z)}{Z}\bigg)\bigg).
\end{align}

\noindent
(ii) Next, we consider the second trace 
\begin{align*}
\rho\tr\bigg(\dpd{}{Z}P\,Q\bigg)=\det(Z-\overline{Z})^{\alpha+\beta-(n+1)}\overline{\psi}(Z)\tr\bigg(\dpd{}{Z}\Lambda_{\beta}\varphi(Z)(Z-\overline{Z})\bigg)
\end{align*}
in equation~\eqref{rhopq0}, which, again through rearrangement of the terms inside the trace by a cyclical permutation, takes the form
\begin{align}
\label{rhopq2}
\rho\tr\bigg(\dpd{}{Z}P\,Q\bigg)=\det(Z-\overline{Z})^{\alpha+\beta-(n+1)}\tr\bigg((Z-\overline{Z})\dpd{}{Z}\Lambda_{\beta}\varphi(Z)\overline{\psi}(Z)\bigg).
\end{align}

\noindent
(iii) Finally, we consider the third trace 
\begin{align*}
\rho\tr\bigg(P^{t}\dpd{}{Z}Q\bigg)=\det(Z-\overline{Z})^{\alpha+\beta-(n+1)}\overline{\psi}(Z)\tr\bigg(\big(\Lambda_{\beta}\varphi(Z)\big)^{t}\bigg(\dpd{}{Z}(Z-\overline
{Z})\bigg)\mathbbm{1}_{n}\bigg)
\end{align*}
in equation~\eqref{rhopq0}. By the first operator identity in Corollary~\ref{operator derivative identity}, we have the matrix identity 
\begin{align*}
\bigg(\dpd{}{Z}(Z-\overline{Z})\bigg)\mathbbm{1}_{n}=\bigg((Z-\overline{Z})\dpd{}{Z}\bigg)^{t}\mathbbm{1}_{n}+\frac{1}{2}(n+1)\mathbbm{1}_{n}=\frac{1}{2}(n+1)
\mathbbm{1}_{n},
\end{align*}
which gives, upon rearrangement of the scalar quantities, the identity
\begin{align}
\label{rhopq3}
\rho\tr\bigg(P^{t}\dpd{}{Z}Q\bigg)=\det(Z-\overline{Z})^{\alpha+\beta-(n+1)}\frac{1}{2}(n+1)\tr\big(\Lambda_{\beta}\varphi(Z)\overline{\psi}(Z)\big).
\end{align}

\noindent
Now, adding up equations \eqref{rhopq1}--\eqref{rhopq3}, it follows from equation~\eqref{rhopq0} that
\begin{align*}
&\sum\limits_{j,k,l=1}^{n}\bigg(\dpd{}{Z}\bigg)_{l,j}(\rho\,p_{j,k}\,q_{k,l})=\det(Z-\overline{Z})^{\alpha+\beta-(n+1)}\bigg((\alpha+\beta-(n+1)/2)\tr\big(\Lambda_{\beta}
\varphi(Z)\overline{\psi}(Z)\big) \\
&\hspace*{35mm}+\tr\bigg(\Lambda_{\beta}\varphi(Z)(Z-\overline{Z})\dpd{\overline{\psi}(Z)}{Z}\bigg)+\tr\bigg((Z-\overline{Z})\dpd{}{Z}\Lambda_{\beta}\varphi(Z)
\overline{\psi}(Z)\bigg)\bigg).
\end{align*}
Rearranging terms on the right-hand side of the last expression, leads to
\begin{align*}
\sum\limits_{j,k,l=1}^{n}\bigg(\dpd{}{Z}\bigg)_{l,j}(\rho\,p_{j,k}\,q_{k,l})&=\det(Z-\overline{Z})^{\alpha+\beta-(n+1)}\bigg(\tr\bigg(\Lambda_{\beta}\varphi(Z)\bigg
((Z-\overline{Z})\dpd{}{Z}+\beta\mathbbm{1}_{n}\bigg)\overline{\psi}(Z)\bigg) \\
&\quad+\tr\bigg((Z-\overline{Z})\dpd{}{Z}+(\alpha-(n+1)/2)\mathbbm{1}_{n}\bigg)\Lambda_{\beta}\varphi(Z)\overline{\psi}(Z)\bigg).
\end{align*}
Identifying the operator $(Z-\overline{Z})\partial/\partial Z+\beta\mathbbm{1}_{n}$ on the right-hand side of the above equation as $-\overline{\Lambda}_{\beta}$ 
and the operator $(Z-\overline{Z})\partial/\partial Z+(\alpha-(n+1)/2)\mathbbm{1}_{n}$ as $K_{\alpha-(n+1)/2}$, we can rewrite the right-hand side of the above 
equation as
\begin{align*}
\det(Z-\overline{Z})^{\alpha+\beta-(n+1)}\Big(-\tr\big(\Lambda_{\beta}\varphi(Z)\,\overline{\Lambda}_{\beta}\overline{\psi}(Z)\big)+\tr\big(K_{\alpha-(n+1)/2}\,\Lambda_
{\beta}\varphi(Z)\overline{\psi}(Z)\big)\Big),
\end{align*}
which, by definition of $\widetilde{\Omega}_{\alpha,\beta}$, is equal to
\begin{align*}
\det(Z-\overline{Z})^{\alpha+\beta-(n+1)}\Big(\tr\big(\widetilde{\Omega}_{\alpha,\beta}-\beta(\alpha-(n+1)/2)\mathbbm{1}_{n}\big)\varphi(Z)\overline{\psi}(Z)-\tr\big
(\Lambda_{\beta}\varphi(Z)\overline{\Lambda}_{\beta}\overline{\psi}(Z)\big)\Big).
\end{align*}
In total, we get
\begin{align*}
\sum\limits_{j,k,l=1}^{n}\bigg(\dpd{}{Z}\bigg)_{l,j}(\rho\,p_{j,k}\,q_{k,l})&=\det(Z-\overline{Z})^{\alpha+\beta-(n+1)}\Big(-\Delta_{\alpha,\beta}\,\varphi(Z)\overline{\psi}
(Z)-\tr\big(\Lambda_{\beta}\varphi(Z)\overline{\Lambda}_{\beta}\overline{\psi}(Z)\big) \\
&\quad-n\beta(\alpha-(n+1)/2)\,\varphi(Z)\overline{\psi}(Z)\Big).
\end{align*}
Thus, substituting $\sum\limits_{j,k,l=1}^{n}(\partial/\partial Z)_{l,j}(\rho\,p_{j,k}\,q_{k,l})$ back into equation \eqref{dabba}, we arrive at
\begin{align*}
\mathrm{d}\omega(Z)&=\det(Z-\overline{Z})^{\alpha+\beta}\Big(-\Delta_{\alpha,\beta}\,\varphi(Z)\overline{\psi}(Z)-\tr\big(\Lambda_{\beta}\varphi(Z)\overline
{\Lambda}_{\beta}\overline{\psi}(Z)\big) \\[1mm]
&\quad-n\beta(\alpha-(n+1)/2)\,\varphi(Z)\overline{\psi}(Z)\Big)\frac{[\mathrm{d}Z]\wedge[\mathrm{d}\overline{Z}]}{\det(Z-\overline{Z})^{n+1}}.
\end{align*}
Now, noting that the volume form 
\begin{align*}
\det(Z-\overline{Z})^{\alpha+\beta}\frac{[\mathrm{d}Z]\wedge[\mathrm{d}\overline{Z}]}{\det(Z-\overline{Z})^{n+1}}
\end{align*}
is just a constant multiple of $\det(Y)^{\alpha+\beta}\mathrm{d}\mu_{n}(Z)$, it follows readily from the vanishing result~\eqref{int zero} that
\begin{align*}
\langle-\Delta_{\alpha,\beta}\varphi,\psi\rangle=\int\limits_{\Gamma\backslash\mathbb{H}_{n}}\det(Y)^{\alpha+\beta}\tr\big(\Lambda_{\beta}\varphi(Z)\overline{\Lambda}_
{\beta}\overline{\psi}(Z)\big)\,\mathrm{d}{\mu}_{n}(Z)+n\beta(\alpha-(n+1)/2)\langle\varphi,\psi\rangle, 
\end{align*}
which is the claimed formula.

Using the latter formula, we compute
\begin{align*}
\langle\varphi,-\Delta_{\alpha,\beta}\psi\rangle&=\overline{\langle-\Delta_{\alpha,\beta}\psi,\varphi\rangle} \\
&=\int\limits_{\Gamma\backslash\mathbb{H}_{n}}\det(Y)^{\alpha+\beta}\,\overline{\tr\big(\Lambda_{\beta}\psi(Z)\overline{\Lambda}_{\beta}\overline{\varphi}(Z)\big)}\,
\mathrm{d}\mu_{n}( Z)+n\beta(\alpha-(n+1)/2)\overline{\langle\psi,\varphi\rangle} \\
&=\int\limits_{\Gamma\backslash\mathbb{H}_{n}}\det(Y)^{\alpha+\beta}\tr\big(\Lambda_{\beta}\varphi(Z)\overline{\Lambda}_{\beta}\overline{\psi}(Z)\big)\,\mathrm{d}
\mu_{n}(Z)+n\beta(\alpha-(n+1)/2)\langle\varphi,\psi\rangle \\[1mm]
&=\langle-\Delta_{\alpha,\beta}\varphi,\psi\rangle,
\end{align*}
which proves the claimed symmetry of the Siegel--Maa\ss{} Laplacian $\Delta_{\alpha,\beta}$.
\end{proof}

\begin{corollary}
\label{characterization}
Let $\Gamma\subset\mathrm{Sp}_{n}(\mathbb{R})$ be a subgroup commensurable with $\mathrm{Sp}_{n}(\mathbb{Z})$ and let $\varphi\in\mathcal{H}_{\alpha,
\beta}^{n}(\Gamma)$ be a Siegel--Maa\ss{} form of weight $(\alpha,\beta)$ and degree $n$ for $\Gamma$. Then, if $\varphi$ is an eigenform of $\Delta_{\alpha,
\beta}$ with eigenvalue $\lambda$, we have $\lambda\in\mathbb{R}$ and $\lambda\geq n\beta(\alpha-(n+1)/2)$.

Furthermore, $\varphi$ has eigenvalue $\lambda=\beta(\alpha-(n+1)/2)$ if and only if $\varphi(Z)=\det(Y)^{-\beta}f(Z)$, where $f\colon\mathbb{H}_{n}\longrightarrow
\mathbb{C}$ is a holomorphic function satisfying
\begin{align*}
f(\gamma Z)=\det(CZ+D)^{\alpha-\beta} f(Z)
\end{align*}
for all $\gamma=\big(\begin{smallmatrix}A&B\\C&D\end{smallmatrix}\big)\in\Gamma$. Moreover, if $\beta<0$, then $f$ is a Siegel cusp form of weight $\alpha-
\beta$ and degree $n$ for~$\Gamma$.
\end{corollary}
\begin{proof}
Since $\varphi\in\mathcal{H}_{\alpha,\beta}^{n}(\Gamma)$ is an eigenform of $\Delta_{\alpha,\beta}$ with eigenvalue $\lambda$, i.e., we have $(\Delta_{\alpha,\beta}+
\lambda\,\mathrm{id})\varphi=0$, we compute using Theorem~\ref{main theorem}
\begin{align*}
\lambda\langle\varphi,\varphi\rangle&=\langle-\Delta_{\alpha,\beta}\varphi,\varphi\rangle \\
&=\int\limits_{\Gamma\backslash\mathbb{H}_{n}}\det(Y)^{\alpha+\beta}\tr\big(\vert\Lambda_{\beta}\varphi(Z)\vert^{2}\big)\,\mathrm{d}\mu_{n}(Z)+n\beta(\alpha-(n+1)/2)
\langle\varphi,\varphi\rangle.
\end{align*}
This immediately implies that $\lambda\in\mathbb{R}$. Furthermore, since $\tr(\vert\Lambda_{\beta}\varphi(Z)\vert^{2})\geq 0$, we conclude that 
\begin{align*}
\lambda\geq n\beta(\alpha-(n+1)/2).
\end{align*}

To prove the second part of the corollary,  we observe that the above equation shows that the equality $\lambda=n\beta(\alpha-(n+1)/2)$ is equivalent to
\begin{align*}
\int\limits_{\Gamma\backslash\mathbb{H}_{n}}\det(Y)^{\alpha+\beta}\tr\big(\vert\Lambda_{\beta}\varphi(Z)\vert^{2}\big)\,\mathrm{d}\mu_{n}(Z)=0.
\end{align*}
Since $\tr(\vert\Lambda_{\beta}\varphi(Z)\vert^{2})\geq 0$, the above integral vanishes if and only if $\tr(\vert\Lambda_{\beta}\varphi(Z)\vert^{2})=0$. Now, as  
the matrix 
\begin{align*}
\Lambda_{\beta}\varphi(Z)=(Z-\overline{Z})\dpd{\varphi}{\overline{Z}}-\beta\varphi(Z)\mathbbm{1}_{n}
\end{align*}
is similar to the complex symmetric matrix 
\begin{align*}
S(Z):=2i\,Y^{1/2}\dpd{\varphi}{\overline{Z}}Y^{1/2}-\beta\varphi(Z)\mathbbm{1}_{n},
\end{align*}
as we have the relation $\Lambda_{\beta}\varphi(Z)=Y^{1/2}S(Z)Y^{-1/2}$, the matrix $\vert\Lambda_{\beta}\varphi(Z)\vert^{2}$ becomes similar to the positive
semidefinite hermitian matrix $S(Z)\overline{S}(Z)$, which is diagonalizable with non-negative real eigenvalues. Therefore, the condition $\tr(S(Z)\overline{S}(Z))=
\tr(\vert\Lambda_{\beta}\varphi(Z)\vert^{2})=0$ is equivalent to the vanishing of all the eigenvalues of $S(Z)\overline{S}(Z)$, which is equivalent to the vanishing 
of $S(Z)$ and hence of $\Lambda_{\beta}\varphi(Z)$. All in all, this proves that the equality $\lambda=n\beta(\alpha-(n+1)/2)$ is equivalent to the vanishing condition
 $\Lambda_{\beta}\varphi=0$.

Continuing, we now set $f(Z):=\det(Y)^{\beta}\varphi(Z)$, and compute
\begin{align*}
\dpd{f}{\overline{Z}}=\beta\det(Y)^{\beta-1}\dpd{\det(Y)}{\overline{Z}}\varphi(Z)+\det(Y)^{\beta}\dpd{\varphi}{\overline{Z}}. 
\end{align*}
Since we have
\begin{align*}
\dpd{\det(Y)}{\overline{Z}}=\frac{1}{2}\bigg(\dpd{}{X}+i\dpd{}{Y}\bigg)\det(Y)=\frac{i}{2}\dpd{\det(Y)}{Y}=\frac{i}{2}\det(Y)Y^{-1},
\end{align*}
the above equality becomes
\begin{align*}
\dpd{f}{\overline{Z}}&=\frac{i\beta}{2}\det(Y)^{\beta}Y^{-1}\varphi(Z)+\det(Y)^{\beta}\dpd{\varphi}{\overline{Z}} \\
&=-\frac{i}{2}\det(Y)^{\beta}Y^{-1}\bigg(-\beta\varphi(Z)\mathbbm{1}_{n}+2iY\dpd{\varphi}{\overline{Z}}\bigg) \\
&=-\frac{i}{2}\det(Y)^{\beta}Y^{-1}\bigg((Z-\overline{Z})\dpd{\varphi}{\overline{Z}}-\beta\varphi(Z)\mathbbm{1}_{n}\bigg) \\
&=-\frac{i}{2}\det(Y)^{\beta}Y^{-1}\Lambda_{\beta}\varphi(Z).
\end{align*}
In total, this shows that $\partial f/\partial\overline{Z}=0$, i.e., the function $f$ is holomorphic, if and only if $\Lambda_{\beta}\varphi(Z)=0$, which, by the previous
argument, is equivalent to $\varphi\in\mathcal{H}_{\alpha,\beta}^{n}(\Gamma)$ being a Siegel--Maa\ss{} form with eigenvalue $\lambda=\beta(\alpha-(n+1)/2)$.
 
Furthermore, as the function $\varphi\in\mathcal{H}_{\alpha,\beta}^{n}(\Gamma)$ has the transformation behaviour
\begin{align*}
\varphi(\gamma Z)=\det(CZ+D)^{\alpha}\det(C\overline{Z}+D)^{\beta}\varphi(Z)
\end{align*}
for all $\gamma=\big(\begin{smallmatrix}A&B\\C&D\end{smallmatrix}\big)\in\Gamma$, the function $f(Z)=\det(Y)^{\beta}\varphi(Z)=\det(\im(Z))^{\beta}\varphi(Z)$ 
has the transformation behaviour
\begin{align*}
f(\gamma Z)&=\det(\im(\gamma Z))^{\beta}\varphi(\gamma Z) \\
&=\bigg(\frac{\det(\im(Z))}{\vert\det(CZ+D)\vert^{2}}\bigg)^{\beta}\det(CZ+D)^{\alpha}\det(C\overline{Z}+D)^{\beta}\varphi(Z) \\[1mm]
&=\det(CZ+D)^{\alpha-\beta}\det(\im(Z))^{\beta}\varphi(Z) \\[2mm]
&=\det(CZ+D)^{\alpha-\beta}f(Z),
\end{align*}
as claimed. 

Finally, letting $\gamma_{j}=\big(\begin{smallmatrix}A_{j}&B_{j}\\C_{j}&D_{j}\end{smallmatrix}\big)\in\mathrm{Sp}_{n}(\mathbb{Z})$ ($j=1,\ldots,h$) be a set of 
representatives for the left cosets of $\Gamma\cap\mathrm{Sp}_{n}(\mathbb{Z})$ in $\mathrm{Sp}_{n}(\mathbb{Z})$, Remark~\ref{boundedness} shows that
given $Y_{0}\in\mathrm{Sym}_{n}(\mathbb{R})$ with $Y_{0}>0$, the quantities
\begin{align*}
\vert\det(C_{j}Z+D_{j})^{-\alpha}\det(C_{j}\overline{Z}+D_{j})^{-\beta}\varphi(\gamma_{j} Z)\vert
\end{align*}
have to be bounded in the region $\{Z=X+iY\in\mathbb{H}_{n}\,\vert\,Y\geq Y_{0}\}$. Therefore, if $\beta<0$, this implies that given $Y_{0}\in\mathrm{Sym}_{n}
(\mathbb{R})$ with $Y_{0}\gg 0$, the quantities
\begin{align*}
\vert\det(C_{j}Z+D_{j})^{-(\alpha-\beta)}f(\gamma_{j} Z)\vert=\vert\det(C_{j}Z+D_{j})^{-\alpha}\det(C_{j}\overline{Z}+D_{j})^{-\beta}\det(\im(\gamma_{j} Z))^{\beta}
\varphi(\gamma_{j} Z)\vert
\end{align*}
will become arbitrarily small in the region $\{Z=X+iY\in\mathbb{H}_{n}\,\vert\,Y\geq Y_{0}\}$. In other words, $f$ is indeed a Siegel cusp form of weight $\alpha-
\beta$ and degree $n$ for $\Gamma$.

With all this, the proof of the corollary is complete.
\end{proof}

\begin{remark}
For $\Gamma\subset\mathrm{Sp}_{n}(\mathbb{R})$ a subgroup commensurable with $\mathrm{Sp}_{n}(\mathbb{Z})$ and $\alpha=k/2$, $\beta=-k/2$ with 
$k\in\mathbb{N}_{>0}$, we denote the Hilbert space $\mathcal{H}_{\alpha,\beta}^{n}(\Gamma)$ simply by $\mathcal{H}_{k}^{n}(\Gamma)$. Similarly, we write 
for the operator $\Omega_{\alpha,\beta}$ simply $\Omega_{k}$, which becomes
\begin{align*}
\Omega_{k}&=(Z-\overline{Z})\bigg((Z-\overline{Z})\dpd{}{\overline{Z}}\bigg)^{t}\dpd{}{Z}+\frac{k}{2}(Z-\overline{Z})\dpd{}{\overline{Z}}+\frac{k}{2}(Z-\overline{Z})
\dpd{}{Z} \\
&=-Y\bigg(\bigg(Y\dpd{}{X}\bigg)^{t}\dpd{}{X}+\bigg(Y\dpd{}{Y}\bigg)^{t}\dpd{}{Y}\bigg)+ikY\dpd{}{X}.
\end{align*}
Finally, we write for the operator $\Delta_{\alpha,\beta}$ simply $\Delta_{k}$ and call it the Siegel--Maa\ss{} Laplacian of weight $k$; it is given as
\begin{align*}
\Delta_{k}=\tr\bigg(Y\bigg(\bigg(Y\dpd{}{X}\bigg)^{t}\dpd{}{X}+\bigg(Y\dpd{}{Y}\bigg)^{t}\dpd{}{Y}\bigg)-ikY\dpd{}{X}\bigg).
\end{align*}
We note that the transformation behaviour of a Siegel--Maa\ss{} form $\varphi$ of weight $k$ and degree $n$ for $\Gamma$ takes the form
\begin{align*}
\varphi(\gamma Z)=\bigg(\frac{\det(CZ+D)}{\det(C\overline{Z}+D)}\bigg)^{k/2}\varphi(Z),
\end{align*}
where $\gamma=\big(\begin{smallmatrix}A&B\\C&D\end{smallmatrix}\big)\in\Gamma$.
\end{remark} 

In the last corollary, we summarize the main results about Siegel--Maa\ss{} forms of weight $k$ and degree~$n$ for $\Gamma$.

\begin{corollary}
\label{kernel connection}
Let $\Gamma\subset\mathrm{Sp}_{n}(\mathbb{R})$ be a subgroup commensurable with $\mathrm{Sp}_{n}(\mathbb{Z})$ and let $\varphi\in\mathcal{H}_{k}^{n}
(\Gamma)$ be a Siegel--Maa\ss{} form of weight $k$ and degree $n$ for $\Gamma$. Then, if $\varphi$ is an eigenform of $\Delta_{k}$ with eigenvalue $\lambda$, 
we have $\lambda\in\mathbb{R}$ and 
\begin{align*}
\lambda\geq \frac{nk}{4}(n-k+1),
\end{align*}
with equality attained if and only if the function $\varphi$ is of the form $\varphi(Z)=\det(Y)^{k/2}f(Z)$ for some Siegel cusp form $f\in\mathcal{S}_{k}^{n}(\Gamma)$ 
of weight $k$ and degree $n$ for $\Gamma$. In other words, there is an isomorphism
\begin{align*}
\mathcal{S}_{k}^{n}(\Gamma)\cong\ker\bigg(\Delta_{k}+\frac{nk}{4}(n-k+1)\mathrm{id}\bigg)
\end{align*}
of $\mathbb{C}$-vector spaces, induced by the assignment $f\mapsto\det(Y)^{k/2}f$.  
\end{corollary}
\begin{proof}
The proof is an immediate consequence of Corollray~\ref{characterization} by setting $\alpha=k/2$ and $\beta=-k/2$.
\end{proof}

\newpage

\bibliographystyle{amsplain}
\bibliography{bibliogen}

\end{document}